\documentclass[12pt]{amsart}
\usepackage{geometry} 
\newcounter{dummy} \numberwithin{dummy}{section}
\newtheorem{theorem}[dummy]{Theorem}
\newtheorem{definition}[dummy]{Definition}
\newtheorem{proposition}[dummy]{Proposition}
\newtheorem{lemma}[dummy]{Lemma}
\usepackage{graphicx}
\usepackage{amssymb}
\usepackage{epstopdf}
\usepackage{lscape}
\usepackage{indentfirst}
\usepackage{latexsym}
\usepackage{multirow}
\usepackage{tabls}
\usepackage{amsthm}
\usepackage{amsthm,amsmath}
\usepackage{graphicx}
\usepackage{lscape}
\usepackage{indentfirst}
\usepackage{latexsym}
\usepackage{multirow}
\usepackage{tabls}
\usepackage{wrapfig}
\usepackage[all]{xy}
\usepackage[numbers]{natbib}

\title[Class Numbers of Real Cyclotomic Fields of Conductor $pq$]{on the class numbers of real cyclotomic fields of conductor $pq$}
\author{eleni agathocleous}

\begin{document}
\begin{center}
\maketitle{\textbf{Abstract}}
\end{center}

\tiny
The class numbers $h^{+}$ of the real cyclotomic fields are very
hard to compute. Methods based on discriminant bounds become useless
as the conductor of the field grows and methods employing Leopoldt's decomposition of the class number become hard to use when the 
field extension is not cyclic of prime power.  This is why 
other methods have been developed which approach the problem from different
angles. In this paper we extend one of these methods that was
designed for real cyclotomic fields of prime conductor, and we make it applicable to real
cyclotomic fields of conductor equal to the product of two distinct odd primes. The main advantage of this method is that it does not
exclude the primes dividing the order of the Galois group, in contrast to other methods. We applied our algorithm to real cyclotomic fields of conductor $<$ 2000 and we calculated the full order of the $l$-part of $h^{+}$ for all odd primes $l$ $<$ 10000.
 \\ 
 
 \textbf{AMS Mathematics Subject Classification:} 11Y40 \\

\textbf{Keywords:} computation of class numbers, real cyclotomic fields, non-prime conductor \\

\normalsize

\section{Introduction}
\label{Intro}
Let $\mathbb{Q}(\zeta_{m})$ be the cyclotomic field of conductor $m$ and
denote by $C$ its ideal class group and by $h = |C|$ its class
number. In the same way let $C^{+}$ and $h^{+}$ denote the ideal
class group and class number of the maximal real subfield
$\mathbb{Q}(\zeta_{m})^{+}$. The natural map $C^{+}\ \longrightarrow \ C$ is
an injection \cite[Theorem 4.14]{W1} and we have the well known
result $h$ = $h^{+} \ h^{-}$. The relative class number $h^{-}$ is
easy to compute as there is an explicit and easily computable
formula for its order \cite[Theorem 4.17]{W1}. The number $h^{+}$ however is extremely hard to compute. The class
number formula is not so useful as it requires that the units of
$\mathbb{Q}(\zeta_{m})^{+}$ to be known. Methods that use the classical
Minkowski bound become useless as $m$ grows, and other methods based
on Odlyzko's discriminant bounds (see \cite{Od1} and \cite{Od2}) are
only applicable to fields with small conductor. Masley in
\cite{Masley} computed the class numbers for real abelian fields of
conductor $\leq$ 100 and Van der Linden in \cite{vdLinden} was able
to calculate the class number of a large collection of real abelian
fields of conductor $\leq$ 200. For fields of larger conductor
however, the above methods can not be effective. 
Many of the other methods that were developed, employ the well known Leopoldt's
decomposition of the class number $h^{+}$ of a real abelian field
$K$, see \cite{L1}, which derives from his decomposition of the
cyclotomic units into the product of the cyclotomic units of all
cyclic subfields $K_{\chi}$ of $K$. More specifically, we have that
$h^{+} =$ $Q$ $\prod_{\chi}h_{\chi}$, where the product runs over all
non-trivial characters $\chi$ irreducible over the rationals, each
`class number' $h_{\chi}$ is the index of the cyclotomic units of
$K_{\chi}$ in its full group of units $E_{\chi}$ and $Q$ is some
value which equals 1 in the case where the extension $K/\mathbb{Q}$ is cyclic
of prime order, but which is very hard to compute in the general
case.

A different method is introduced by Schoof in \cite{Schoof1} and
is designed for real cyclotomic fields of prime conductor. Schoof
developed an algorithm that computes the order of the module
$B$=Units/(Cyclotomic Units), which is precisely equal to $h^{+}$ in
his case, where the conductor of the field is a prime number. One of the great advantages of his method is that it does
not exclude the primes dividing the order of the extension, as opposed to other methods.

In this paper we extend Schoof's method to real
cyclotomic fields of conductor equal to the product of two distinct odd primes. We apply our algorithm to real cyclotomic fields of conductor $<$ 2000 and we calculate the $l$-part of $h^{+}$ for all odd primes $l$ $<$ 10000.

\section{Our Cyclotomic Unit $\eta$}
\label{eta}
In our case, where the conductor of the field is not a prime number, the group of cyclotomic units has a complicated structure. We therefore work with a cyclic subgroup, yet of finite index in the group of units, which we present below.

Let $p$ and $q$ be distinct odd primes. From now on, $G$ will denote the Galois group $Gal(\mathbb{Q}(\zeta_{pq})^{+}/\mathbb{Q})$, $E$ will denote 
the group of units of the real cyclotomic field $\mathbb{Q}(\zeta_{pq})^{+}$
and $O$ its ring of integers. Without loss of generality we will
always assume that $p$ $<$ $q$. Choose and fix $g$ and $h$,
primitive roots modulo $p$ and $q$ respectively. Denote by
$\eta_{(g,h)}$ the following real unit of $\mathbb{Q}(\zeta_{pq})^{+}$:
$$\eta_{(g,h)} = \zeta_{pq}^{-(p+q)}(1-{\zeta_{pq}}^{p+q})^{2}
\frac{\zeta_{p}^{-g/2}}{\zeta_{p}^{-1/2}} \ \frac{(1 -
\zeta_{p}^{g})}{(1 - \zeta_{p})} \
\frac{\zeta_{q}^{-h/2}}{\zeta_{q}^{-1/2}} \ \frac{(1 -
\zeta_{q}^{h})}{(1 - \zeta_{q})}$$ and by $H_{(g,h)}$ the group $\pm
\eta_{(g,h)}^{\mathbb{Z}[G]}$. We will omit the subscripts and just
write $H$ and $\eta$ since we will let $\eta_{\alpha}$ denote the result of the action of the element $\alpha$ $\in$ $G$ on $\eta$. With this
notation in mind, we are ready to prove a statement about the
regulator of the units $\{\eta_{\alpha}\}_{\alpha \in G}$.

\begin{proposition} \label{prop1} Let $E$ be the group of units of $\mathbb{Q}(\zeta_{pq})^{+}$ and $H$ = $\pm \eta^{\mathbb{Z}[G]}$
= $\pm \eta_{(g,h)}^{\mathbb{Z}[G]}$ as above, where $g$ and $h$ are any two
fixed primitive roots modulo $p$ and $q$ respectively. The index
[$E$:$H$] is always finite and it equals: $$[E:H] =\frac{2^{|G|-1}
h^{+}}{|G|} \cdot$$ $$\prod_{\chi=\chi_{p}\neq1} \frac{1}{2}
\big{[}2(\chi(q)^{-1}-1)+(\chi(g^{-1})-1)(q-1)\big{]}\cdot$$
$$\prod_{\chi=\chi_{q}\neq1}\frac{1}{2}
\big{[}2(\chi(p)^{-1}-1)+(\chi(h^{-1})-1)(p-1)\big{]}$$ where the
characters $\chi$ in the first product are the even characters
$\chi_{p}$ of conductor $p$ and those in the second product are the
even characters $\chi_{q}$ of conductor $q$.
\end{proposition}

\emph{Proof}: Define $f$ by $f(\alpha)$ = $\log|\eta_{\alpha}|$. We
see that $\sum_{\alpha}f(\alpha)$ =
$\log|\prod_{\alpha}\eta_{\alpha}|$ = 0. Denote by $\chi$ an even
Dirichlet character and note that for any root of unity $\zeta$ we
have that $\log|\zeta^{i}(1-\zeta^{j})| = \log|1-\zeta^{j}|$, where
$i$ and $j$ are arbitrary. The regulator $R$ of the units
$\eta_{\alpha}$ is: $$R = R(\{\eta_{\alpha}\}) = \pm \det(\log|\eta_{\alpha\beta}|)_{\alpha,\beta\neq 1} = \pm \det(f(\alpha\beta))_{\alpha,\beta\neq 1}$$
$$=\pm \det(f(\beta\alpha^{-1}))_{\alpha,\beta\neq 1} \phantom{xy} (by \ rearranging \ the \ rows)$$
$$=\pm \frac{1}{|G|}\prod_{\chi \neq 1}\sum_{\beta\in G} \chi(\beta)f(\beta) \phantom{xy} \text{by \citep[Lemma 5.26(c)]{W1}})$$
$$=\pm \frac{1}{|G|}\prod_{\chi \neq 1} \frac{1}{2} \sum_{^{1 \leq \beta \leq pq}_{(\beta,pq)=1}} \chi(\beta)\big{[}\log|1 - \zeta_{pq}^{\beta(p+q)}|^{2} + \log|\frac{1 - \zeta_{p}^{g\beta}}{1 - \zeta_{p}^{\beta}}| + \log|\frac{1 -
\zeta_{q}^{h\beta}}{1 - \zeta_{q}^{\beta}}|\big{]} $$

For a character $\chi$ and the second summand we have $$\sum_{^{1
\leq \beta \leq pq}_{(\beta,pq)=1}} \chi(\beta) \big{[}\log|1 \ - \
\zeta_{p}^{g\beta}|\ - \\log|1 - \zeta_{p}^{\beta}|\big{]} = 0$$
for $f_{\chi} = pq$ by \text{\citep[Lemma 8.4]{W1}} and for $f_{\chi} = q$ by \text{ \citep[Lemma 8.4 and 8.5]{W1}}.

For $f_{\chi} = p$ and by applying \text{\citep[Lemma 8.4 and 8.5]{W1}}, the above sum equals

$$=\chi(g)^{-1}\ \sum_{^{g\beta (mod \ pq)}_{(\beta,pq)=1}} \chi(g
\beta)\log|1 - \zeta_{p}^{g\beta}| \ - \ \sum_{^{1 \leq \beta \leq
pq}_{(\beta,pq)=1}} \chi(\beta) \log|1\ -\ \zeta_{p}^{\beta}|$$
$$=(\chi(g^{-1})-1) \sum_{^{1 \leq \alpha \leq pq}_{(\alpha,pq)=1}}
\chi(\alpha) \log|1 \ - \ \zeta_{p}^{\alpha}|.$$ 

To sum up, the second summand gives:

\[ = \left\{\begin{array}{r@{\quad,\quad}l} 0 & if f_{\chi} = pq \\(\chi(g^{-1})-1)(q-1)\sum_{1 \leq \alpha \leq p} \chi(\alpha)
\log|1 \ -\ \zeta_{p}^{\alpha}| & if \ f_{\chi} = p \\
0 & if \ f_{\chi} = q
\end{array}\right.
\]

Similarly, the third summand equals:

\[ = \left\{\begin{array}{r@{\quad,\quad}l} 0 & if f_{\chi} = pq \\0 & if \ f_{\chi} = p \\
(\chi(h^{-1})-1)(p-1)\sum_{1 \leq \alpha \leq q} \chi(\alpha) \log|1
\ - \ \zeta_{q}^{\alpha}| & if \ f_{\chi} = q
\end{array}\right.
\]

It is necessary to state here that the summands of the form $0 \cdot \log \ 0$ are treated as $0$ in formulae as the above.

Putting all three together and denoting by $\chi_{pq}$, $\chi_{p}$
and $\chi_{q}$ the characters of conductor $pq$, $p$ and $q$
respectively, we have that
$$R \ = \ \pm \frac{1}{|G|} \prod_{\chi=\chi_{pq} \neq 1} \frac{1}{2} \cdot 2 \ \chi(p+q)^{-1} \sum_{1 \leq \beta \leq pq} \chi(\beta) \ \log|1 -
\zeta_{pq}^{\beta}| \cdot$$
$$\prod_{\chi=\chi_{p} \neq 1} \ \frac{1}{2} \big{[}2 \chi(q)^{-1} \ (1 - \chi(q))+(\chi(g^{-1})-1)(q-1)\big{]}\sum_{1 \leq \alpha \leq p} \chi(\alpha) \log|1 - \zeta_{p}^{\alpha}|\cdot$$
$$\prod_{\chi=\chi_{q} \neq 1} \ \frac{1}{2} \big{[}2 \chi(p)^{-1} \ (1 - \chi(p))+(\chi(h^{-1})-1)(p-1)\big{]}\sum_{1 \leq \alpha \leq q} \chi(\alpha) \log|1 - \zeta_{q}^{\alpha}|$$

To show that [$E$:$H$] is always finite it suffices to show that the regulator is never zero. Assume it is zero. Then for some character of conductor 
$p$ the sum $$2(\chi(q)^{-1} -1)+(\chi(g^{-1})-1)(q-1)$$ is zero or for some character of conductor $q$ the sum
$$2(\chi(p)^{-1} -1)+(\chi(h^{-1})-1)(p-1)$$ is zero. But
$$2(\chi(q)^{-1} -1)+(\chi(g^{-1})-1)(q-1) = 0 \Leftrightarrow 2\chi(q)^{-1} + (q-1) \chi(g)^{-1} = 2 + (q-1)$$ which never
happens as $\chi(g)^{-1}$ can never equal 1, since $g$ is a primitive root. Similarly for a character of conductor $q$.
Therefore, the regulator is never zero and this completes the proof of Proposition 2.1.\ $\Box$

Denote by $P$ the factor $$\frac{2^{|G|-1}}{|G|}\cdot$$
$$\prod_{\chi=\chi_{p} \neq 1} \big{[}2(\chi(q)^{-1}-1)+(\chi(g^{-1})-1)(q-1)\big{]}$$ $$\prod_{\chi=\chi_{q} \neq 1}\big{[}2(\chi(p)^{-1}-1)+(\chi(h^{-1})-1)(p-1)\big{]}$$ 
which appears in the index [$E$:$H$] in Proposition~\ref{prop1} above. We now have
$$[E:H] = P \cdot h^{+}.$$ One can take advantage of the fact that any choice of primitive roots $g$ and $h$ give a finite index, and for
each field $\mathbb{Q}(\zeta_{pq})^{+}$ one can choose the pair ($g$,$h$) with the property that $P_{(g,h)}$ is divisible by the smallest number of distinct primes. Furthermore, for the primes that appear in this $P_{(g,h)}$ one can can check to see if those primes divide the greatest common divisor of all the $P_{(g,h)}$ for every pair of
primitive roots ($g$,$h$). In the case that a prime $l$ does not divide the greatest common divisor, there is some pair $(g_{0},h_{0})$ for which $l$ does not divide 
$P_{(g_{0},h_{0})}$. We can therefore repeat the first part of our algorithm that we present in Section~\ref{Algorithm}, for this pair $(g_{0},h_{0})$ and for this prime $l$. We continue with the next step of the algorithm for this prime $l$, only if it gives a non-trivial factor. These facts are very useful in the computations, since they narrow down the number of primes that one needs to check to see if they divide $h^{+}$.

\section{The module B= E/H = E/$\pm\eta^{\mathbb{Z}[G]}$}  \label{moduleB} We denote by $B$ be the $\mathbb{Z}[G]$-module $E/H$, where $H = \pm\eta^{\mathbb{Z}[G]}$ as above. From Proposition~\ref{prop1} we have that the order of $B$ is finite and equals the index $\mbox{[E:H]}$. Therefore, by generalizing Schoof, we can calculate the order of $B$ and then multiply by $1/P$ in order to get $h^{+}$, as desired. 

Since $H$ is of finite index in $E$ we have that the map $$\Phi : \mathbb{Z}[G] \rightarrow E$$ $$\phantom{x}  \phantom{x} \alpha \phantom{x} \mapsto \eta^{a}$$ is a homomorphism
whose image $H$ is of finite index and therefore $\mathbb{Z}$-isomorphic to $\mathbb{Z}^{|G|-1}$. We have that $H/\{\pm 1\}$ $\cong$ $\mathbb{Z}[G]/N_{G}$ as $\mathbb{Z}[G]$-modules, where $N_{G}$ is the norm of $G$. Let $M$ $>$ 1 denote a power of a prime $l$. We let $F$ = $\mathbb{Q}(\zeta_{pq})^{+}(\zeta_{2M})$ and $\Delta$ = $Gal(F/\mathbb{Q}(\zeta_{pq})^{+})$. 

The following Lemma appears in \cite{Schoof1}. We prove it again here as we found a few ommisions in the proof. 

\begin{lemma} \label{LemmaKernel} The kernel of the natural map \begin{center}$j : E/E^{M} \rightarrow F^{*}/{F^{*}}^{M}$\end{center} is trivial if $l$ odd and it has order two and is generated by $-1$ if $l$ = 2. \end{lemma}
\emph{Proof:} Fix an embedding $F$ $\subset$ $C$. Then $\mathbb{Q}(\zeta_{pq})^{+}$ identifies with a subfield of $R$. Suppose 0 $<$
$x$ $\in$ $E$ $\subset$ $R$ is in $\operatorname{Ker}j$. Then $x$ = $y^{M}$, some $y$ $\in$ $F^{*}$. Since $\mu_{M} \subset F$ we may assume that $y$
$\in$ R and therefore $conj(y)$=$y$, where $conj$ is complex conjugation in $\Delta$. Since $\Delta$ commutative, $s(y)$=$s(conj(y))$=$conj(s(y))$ $\forall s \in \Delta$, therefore $s(y)$=$\pm y$ $\forall s \in \Delta$, since $y$ and all its conjugates are real $M$-th roots of $x$. If $l \neq 2$ then $M$ is odd. Assume $\exists s \in \Delta$ 
with $s(y)$=$-y$. Then $x$=$s(x)$=$s(y^{M})$=$(s(y))^{M}$=$(-y)^{M}$=$-x$, contradiction. Therefore $\Delta$ fixes $y$ and hence $y$ $\in$ $(\mathbb{Q}(\zeta_{pq})^{+})^{*}$ and
$x$ $\in E^{M}$. Since we took $x$ $>$ 0 we need to check for $-1$ as well. Since $M$ odd, $(-1)^{M}$ = $-1$ therefore $-1$ $\in E^{M}$ as well and in this case $j$ is an
injection. If $l$ = 2 we see that $s(y^{2})$ = $s(y)^{2}$ = $y^{2}$ therefore $y^{2}$ $\in$ ${\mathbb{Q}(\zeta_{pq})^{+}}^{*}$. The quadratic subextensions of $F/\mathbb{Q}(\zeta_{pq})^{+}$ are 
$\mathbb{Q}(\zeta_{pq})(i)$ and $\mathbb{Q}(\zeta_{pq})^{+}(\sqrt{\pm 2}$) and hence $y^{2}$ = $2u^{2}$ or = $\pm u^{2}$, for some $u$ $\in$ $E$. If $y^{2}$ = $2u^{2}$, then 2 = $y^{2}$ $v^{2}$, with $v$ such that $vu$ =1, which can not happen since then (2) = $(v)^{2}$ as ideals but 2 does not ramify in $\mathbb{Q}(\zeta_{pq})^{+}$. So we can only have the second case where
$x$=$y^{M}$=($y^{2}$)$^{2^{(k-1)}}$. For $k \geq 2$ we have $x = u^{2}$ and therefore $x$ $\in$ $E^{M}$. When $k = 1$ we have $x = y^{2} = \pm u^{2}$, but since $x > 0$ we still get $x = u^{2}$ which implies that $x \in E^{M}$. For $-1$, observe that $-1$ = $\zeta_{2M}^{M}$ but $-1$ is not even a square in $\mathbb{Q}(\zeta_{pq})^{+}$ which means 
$\operatorname{Ker}j$ = $\langle -1 \rangle$ of order two in this case. $\Box$

Let $\Omega$ = $Gal(F/\mathbb{Q})$. We have the following exact sequence of galois groups $$0 \rightarrow \Delta \rightarrow \Omega \rightarrow G \rightarrow 0.$$ Let $\Re$ be any prime ideal of $F$ of degree 1, $\rho$ a prime ideal of $\mathbb{Q}(\zeta_{pq})^{+}$ and $r$ a prime number such that $\Re$ $|$ $\rho$ $|$ $r$. We have 
$r\equiv \pm1~({\rm mod}~pq)$ and $r\equiv 1~({\rm mod}~2M)$. We have $|G|$ = $\frac{(p-1)(q-1)}{2}$ and we consider the following diagram:

\footnotesize
$$\xymatrix{ \varepsilon \in E \ar[r]^{f_{1}\phantom{wxyz}} &
\vec{\varepsilon} \in (O/rO)^{*} \ar[d]^{f_{3}}\ar[r]^{f_{2}}
&{(O_{F}/rO_{F})^{*}}^{\Delta}\ar[d]^{f'_{3}}\\
&\mu_{M}(O/rO)&\mu_{M}(O_{F}/rO_{F})^{\Delta}\ar[l]^{f'_{2}}&(\mathbb{Z}/M\mathbb{Z})[\Omega]^{\Delta}\ar[l]^{\phantom{xy}f_{4}}&(\mathbb{Z}/M\mathbb{Z})[G]\ar[l]^{\phantom{xy}f_{5}}}$$  \\
\normalsize

The maps $f_{1}, f_{2}, f'_{2}, f_{3}, f'_{3}, f_{4}$ and $f_{5}$ are defined as in \cite{Schoof1}. 
Let $f_{\Re}$ = $f_{5}^{-1} f_{4}^{-1} f'_{3} f_{2} f_{1}$. In a similar manner as in the proof of Theorem 2.2 in \cite{Schoof1}, it can be shown that the maps $f_{\Re}$ correspond to the Frobenius elements of the primes over $\Re$ in $Gal(F(\sqrt[M]{E})/F)$.
Furthermore, every map in $Hom_{R}(E/\pm E^{M}, R)$ is of the form $f_{\Re}$ for some $\Re \in \ S$, where $S$ denotes the set of unramified prime ideals 
$\Re$ of $\mathbb{Q}(\zeta_{pq})^{+}(\zeta_{2M})$ of degree 1 and $R = (\mathbb{Z}/M\mathbb{Z})[G]$. We can therefore state the following theorem, whose proof we omit as it is very similar to that of Theorem 2.2 in \cite{Schoof1}. 

\begin{theorem} \label{frobenius} Let $l$ and $M$ be as above and let $I$ denote the augmentation ideal of $R$ = $(\mathbb{Z}/M\mathbb{Z})[G]$. We have 
$B[M]^{\perp}$ $\cong$ $I / {\{f_{\Re} (\eta) : \Re \in S\}}$.\end{theorem}

\section{The computations} \label{computations}

\subsection{Reformulating Theorem~\ref{frobenius} in terms of Polynomials} \label{TheoremPolys}
Let $l$ be a fixed odd prime, $M$ $>$ 1 some fixed power of $l$ and $G$ denotes the galois group of $\mathbb{Q}(\zeta_{pq})^{+}$. The group $G$
is of order $(p-1)(q-1)/2$ and we have the isomorphisms 
$$G \cong \Big{(}(\mathbb{Z}/p\mathbb{Z})^{\star}\times(\mathbb{Z}/q\mathbb{Z})^{\star}\Big{)}/\{\pm \mathbf{1}\} \cong$$ 
$$\langle \sigma, \tau : \sigma^{(p-1)} = 1, \tau^{(q-1)} =1, \sigma^{(p-1)/2} \tau^{(q-1)/2} = 1, \sigma \tau = \tau \sigma \rangle$$ where 
$\sigma$ : $\zeta_{p} \mapsto \zeta_{p}^{\gamma}$ and $\tau$ : $\zeta_{q} \mapsto \zeta_{q}^{\delta}$ with $\gamma$ and $\delta$ being fixed
primitive roots modulo $p$ and $q$ respectively. The third relation is that of complex conjugation. The primitive roots $\gamma$ and $\delta$ will be fixed throughout and
will always represent the generators of $(\mathbb{Z}/p\mathbb{Z})^{\times}$ and $(\mathbb{Z}/q\mathbb{Z})^{\times}$ respectively. We see that \begin{center}$\mathbb{Z}[G]$ $\cong$ 
$\mathbb{Z}[x,y]/(x^{p-1} - 1, y^{q-1} - 1, x^{(p-1)/2} y^{(q-1)/2} - 1)$\end{center} via the map that sends $\sigma$ to $x$ and $\tau$ to $y$.
Similarly, $$(\mathbb{Z}/M\mathbb{Z})[G] \cong (\mathbb{Z}/M\mathbb{Z})[x,y]/(x^{p-1} - 1, y^{q-1} - 1, x^{(p-1)/2} y^{(q-1)/2} - 1).$$ Using this notation, the maps $f_{\Re}$ that were introduced in the
previous chapter can now be expressed as polynomials in the variables $x$ and $y$ as follows:
$$f_{\Re}(x,y) = \sum_{1 \leq i \leq p-1} \sum_{1 \leq j \leq (q-1)/2} \log_{r} (\eta_{(i,j)}) \cdot x^{i} \cdot y^{j}$$ where
$$\eta_{(i,j)} = \zeta_{p}^{-\gamma^{i}}\zeta_{q}^{-\delta^{j}}(1-\zeta_{p}^{\gamma^{i}}\zeta_{q}^{\delta^{j}})^{2} \frac{\zeta_{p}^{-g\gamma^{i}/2}}{\zeta_{p}^{-\gamma^{i}/2}} \
\frac{(1 - \zeta_{p}^{g\gamma^{i}})}{(1 - \zeta_{p}^{\gamma^i})} \
\frac{\zeta_{q}^{-h\delta^{j}/2}}{\zeta_{q}^{-\delta^{j}/2}} \
\frac{(1 - \zeta_{q}^{h\delta^{j}})}{(1 - \zeta_{q}^{\delta^{j}})}.$$ Here, $\log_{r}$ denotes the discrete $\log$ which gives $\log_{r}(\eta) = s$ where $s \in \mathbb{Z}/M\mathbb{Z}$ is such
that $\eta^{(r-1)/M}\equiv \zeta_{M}^{s}~({\rm mod}~ \Re)$.

We note here that the second sum in the definition of $f_{\Re}(x,y)$ goes from 1 up to $(q-1)/2$ since we are in the real subfield of $\mathbb{Q}(\zeta_{pq})$.

Given the above, we can now reformulate Theorem~\ref{frobenius} of the previous chapter as follows:

\begin{theorem} \label{frobeniusComp} Let $l$ be a fixed prime and let $M$ $>$ 1 be some fixed power of $l$. Denote by $R$ the ring 
$$(\mathbb{Z}/M\mathbb{Z})[x,y] / (x^{p-1} - 1, y^{q-1} - 1, x^{(p-1)/2} y^{(q-1)/2} - 1)$$ and let $B[M]^{\perp}$ be as in Theorem~\ref{frobenius}. Then 
$$B[M]^{\perp} \cong (x-1,y-1) \Big{/} \{f_{\Re}(x,y) : \Re \ \in S \}$$ where $S = \{the \ degree \ 1 \ prime \ ideals \ of \ \mathbb{Q}(\zeta_{pq})^{+}(\zeta_{2M})\}$.\end{theorem}

\emph{Proof:} From our polynomial description of $\mathbb{Z}[G]$ above, it follows that the augmentation ideal of $(\mathbb{Z}/M\mathbb{Z})[G]$ is $(x-1,y-1)$ .
The result is now immediate from Theorem~\ref{frobenius}. \ $\Box$

\subsection{The Decomposition of the modules $B[M]^{\perp}$} \label{Decomposition}

Let $\tilde{G}$ denote the Galois group of the extension $\mathbb{Q}(\zeta_{pq})/\mathbb{Q}$. We can write $\mathbb{Z}_{l}[\tilde{G}]$ as follows: for the same fixed prime $l$ as above, write 
$p-1 = m_{1} l^{a_{1}}$ and $q-1 = m_{2} l^{a_{2}}$ where $l^{a_{1}} || p-1$ and $l^{a_{2}} || q-1$. Since now $l$ does not divide $m_{1}$ and $m_{2}$, we have
that $$\mathbb{Z}_{l}[\tilde{G}] \cong \mathbb{Z}_{l}[x,y]/(x^{p-1} - 1, y^{q-1} - 1) \cong \mathbb{Z}_{l}[x,y]/((x^{l^{a_{1}}})^{m_{1}} - 1, (y^{l^{a_{2}}})^{m_{2}} -1) \cong$$
$$\prod_{\phi_{x},\phi_{y}}\mathbb{Z}_{l}[x,y]/(\phi_{x}(x^{l^{a_{1}}}),\phi_{y}(y^{l^{a_{2}}})) \cong \prod_{\phi_{x}} \mathbb{Z}_{l}[x]/(\phi_{x}(x^{l^{a_{1}}})) \otimes \prod_{\phi_{y}}\mathbb{Z}_{l}[y]/(\phi_{y}(y^{l^{a_{2}}})) $$
 where the product runs over all irreducible divisors $\phi_{x}$ of $x^{m_{1}}-1$ and $\phi_{y}$ of $y^{m_{2}}-1$. We see that
 $\mathbb{Z}_{l}[x]/(\phi_{x}(x^{l^{a_{1}}}))$ and
 $\mathbb{Z}_{l}[y]/(\phi_{y}(y^{l^{a_{2}}}))$ are complete local
 $\mathbb{Z}_{l}[\tilde{G}]$-algebras with maximal ideals $(l,\phi_{x}(x))$ and
 $(l,\phi_{y}(y))$, respectively, and the orders of their residue fields are $l^{f_{1}}$ and $l^{f_{2}}$, where $f_{1}$ = deg($\phi_{x}(x)$) and $f_{2}$
 = deg($\phi_{y}(y)$). Let $\Delta$ denote the subgroup of $\tilde{G}$ of order prime to $l$. From the decomposition of $\mathbb{Z}_{l}[\tilde{G}]$
above, we can write any finite $\mathbb{Z}_{l}[\tilde{G}]$-module $A$ as a product of its $\phi$-parts
$$A_{\phi_{x},\phi_{y}}=A \otimes_{\mathbb{Z}_{l}[\tilde{G}]}\Big{(} \mathbb{Z}_{l}[x]/(\phi_{x}(x^{l^{a_{1}}})) \otimes \mathbb{Z}_{l}[y]/(\phi_{y}(y^{l^{a_{2}}}))\Big{)}.$$ 
The simple Jordan-H\"older factors of each $A_{\phi_{x},\phi_{y}}$ over $\mathbb{Z}_{l}[\Delta]$ are the same as those over $\mathbb{Z}_{l}[\tilde{G}]$ since we `removed' the powers of $x$ and $y$ dividing the order of $\tilde{G}$.

All of the above about the module $A$ also hold in particular for
$B$, the various $B[M]^{\perp}$ and their $\phi$-parts
$B[M]^{\perp}_{\phi_{x},\phi_{y}}$. Therefore, when we want to find
the Jordan-H\"older factors of $B$ we can start by taking all
combinations of degrees $f_{1}$ and $f_{2}$. Since $x$ and $y$ are
non-zero elements in the corresponding residue fields $\Big{[}
\mathbb{Z}_{l}[x]/(\phi_{x}(x^{l^{a_{1}}}))\Big{]}/(l,\phi_{x}(x))$ and
$\Big{[}\mathbb{Z}_{l}[y]/(\phi_{y}(y^{l^{a_{2}}})\Big{]}/(l,\phi_{y}(y))$,
we must have that the orders of $x$ and $y$ in the ring attached to
$\phi_{x}$ and $\phi_{y}$ must divide $l^{f_{1}}-1$ and
$l^{f_{2}}-1$ respectively. Let $d_{1} =  \operatorname{gcd}(p-1,l^{f_{1}}-1)$
and $d_{2} = \operatorname{gcd}(q-1,l^{f_{2}}-1)$ and let $$R_{d_{1},d_{2}} =
(\mathbb{Z}/M\mathbb{Z})[x,y]/((x^{l^{a_{1}}})^{d_{1}}-1,(y^{l^{a_{2}}})^{d_{2}}-1).$$
Since the rings $R_{d_{1},d_{2}}$ and $R_{\phi_{x},\phi_{y}}$ are
direct summands of $R$, any map from their modules
$B[M]_{d_{1},d_{2}}$ and $B[M]_{\phi_{x},\phi_{y}}$, respectively,
to $R$ will end up in these smaller rings. Therefore we can refer to
$B[M]_{d_{1},d_{2}}^{\perp}$ and $B[M]_{\phi_{x},\phi_{y}}^{\perp}$
as $R_{d_{1},d_{2}}$ and $R_{\phi_{x},\phi_{y}}$ modules,
respectively, where
$$B[M]_{d_{1},d_{2}}^{\perp} \cong I_{d} / \langle
(x^{l^{a_{1}}})^{d_{1}}-1,(y^{l^{a_{2}}})^{d_{2}}-1,cnj,f_{\Re}(x,y)
: \Re \in S\rangle \phantom{wxyz} (i)$$ and similarly for $B[M]_{\phi_{x},\phi_{y}}^{\perp}$. Here $cnj$ denotes the conjugation relation already defined in Section~\ref{TheoremPolys}. Basically $B[M]_{d_{1},d_{2}}^{\perp}$ is the direct sum of the $B[M]_{\phi_{x},\phi_{y}}^{\perp}$'s and we therefore have $|B[M]_{d_{1},d_{2}}^{\perp}|$ = $\prod_{\phi_{x},\phi_{y}} |B[M]_{\phi_{x},\phi_{y}}^{\perp}|$.

Since $(1 \pm cnj)/2$ are
idempotents in $(\mathbb{Z}/M\mathbb{Z})[\tilde{G}]$ for $M$ odd, the conjugation relation in the
ideal
$$J = \langle
(x^{l^{a_{1}}})^{d_{1}}-1,(y^{l^{a_{2}}})^{d_{2}}-1,cnj,f_{\Re}(x,y)
: \Re \in S\rangle \phantom{wxyz} (ii)$$ from $(i)$ above, makes
$B[M]_{d_{1},d_{2}}^{\perp}$ a $(\mathbb{Z}/M\mathbb{Z})[G]$-module as well. 

Note that here,
the polynomials $f_{\Re}$ are restrictions of the Frobenius elements
of Theorem~\ref{frobeniusComp} to this smaller extension determined by the set of
polynomials $(x^{l^{a_{1}}})^{d_{1}}-1$ and
$(y^{l^{a_{2}}})^{d_{2}}-1$. They are therefore of the form
$$f_{\Re}(x,y) = \sum_{1 \leq i \leq d_{1}l^{a_{1}}} \sum_{1 \leq j
\leq d_{2}l^{a_{2}}} \log_{r} (\prod_{^{m\equiv i~({\rm mod}~d_{1}l^{a_{1}})}_{n\equiv j~({\rm mod}~d_{2}l^{a_{2}})}} \eta_{(m,n)})
\cdot x^{i} \cdot y^{j} \phantom{wxyz} (iii).$$

\subsection{Gr\"obner Bases} \label{Groebner}

We make use of Gr\"obner Bases, which we present here
by following \cite{Adams}, in order to handle the appearance of
two variables in our calculations of the ideals $J$
defined in the previous section, and enable ourselves to calculate the order of
the various $B[M]_{d_{1},d_{2}}^{\perp}$. 

As before, $d_{1}=\operatorname{gcd}(p-1,l^{f_{1}}-1)$ and $d_{2}=\operatorname{gcd}(q-1,l^{f_{2}}-1)$, where $f_{1}$ and
$f_{2}$ are the degrees of some irreducible polynomials $\phi_{x}$
and $\phi_{y}$ respectively. Again, let $B[M]_{d_{1},d_{2}}^{\perp}$
be the corresponding $R_{d_{1},d_{2}}$-module and $$J = \langle
(x^{l^{a_{1}}})^{d_{1}}-1,(y^{l^{a_{2}}})^{d_{2}}-1,cnj,f_{\Re}(x,y)
: \Re \in S\rangle$$ the corresponding ideal. All the computations
for the calculation of the Frobenius polynomials were performed in
PARI and the computations for a basis for the ideal $J$ in
MATHEMATICA, which allows the computations of bases for ideals whose
elements are polynomials in more than one variable and their
coefficients are in any ring $(\mathbb{Z}/M\mathbb{Z})$, not necessarily a field.

In this section, $R = (\mathbb{Z}/M\mathbb{Z})[x,y]$ will be our polynomial ring in two
variables $x$ and $y$ with coefficients in the Noetherian ring $(\mathbb{Z}/M\mathbb{Z})$, which makes $R$ Noetherian as well. 
Because of the appearance of more than one
variable in our polynomials, we need to agree on the order of the
variables and also find a way to compare every element. We call a
\emph{power product} an element of the form $x^{a}y^{b}$ with $a,b$
non-negative integers and we denote by $T^{2}$ the set of all power
products our polynomial ring $R$.
Following the definition of \emph{term order} given in \cite{Adams},
we define a total order on $T^{2}$ as follows:
\begin{definition} By a term order on $T^{2}$ we mean a total order
$<$ on $T^{2}$ which satisfies the following conditions:

(i) 1 $<$ $x^{a}y^{b}$ for all $1 \neq x^{a}y^{b} \in T^{2}$

(ii) If $x^{a_{1}}y^{b_{1}} < x^{a_{2}}y^{b_{2}}$ then
$x^{a_{1}}y^{b_{1}}x^{c}y^{d} <
x^{a_{2}}y^{b_{2}}x^{c}y^{d}$ for all $x^{c}y^{d}$ $\in$
$T^{2}$. \end{definition}

The type of \emph{term order} that we use here is the
\emph{lexicographical order} on $T^{2}$ which we define below:

\begin{definition} The lexicographical order on $T^{2}$ with $x > y$ is
defined as:

For $(a_{1},b_{1}), (a_{2},b_{2})$ with $a_{i}, b_{i}$ positive
integers, we define $x^{a_{1}}y^{b_{1}} < x^{a_{2}}y^{b_{2}}$ if and
only if ($a_{1} < a_{2}$ or ($a_{1}=a_{2}$ and $b_{1} < b_{2}$)). We
therefore have \begin{center} $ 1 < y < y^{2} < y^{3} < ... < x < xy
< xy^2 < ... < x^{2} < ...$ \end{center} \end{definition}

Now that we have chosen a term order on our polynomial ring, for
each polynomial \begin{center} $f = c_{1}x^{a_{1}}y^{b_{1}} +
c_{2}x^{a_{2}}y^{b_{2}} + ... + c_{n}x^{a_{n}}y^{b_{n}}$
\end{center} with $c_{i} \neq$ 0 in $(\mathbb{Z}/M\mathbb{Z})$ and $x^{a_{1}}y^{b_{1}}
> x^{a_{2}}y^{b_{2}} > ... > x^{a_{n}}y^{b_{n}}$, we can define:

$\operatorname{lp}(f) = x^{a_{1}}y^{b_{1}}$, the \emph{leading power product} of
$f$,

$\operatorname{lc}(f) = c_{1}$, the \emph{leading coefficient} of $f$,

$\operatorname{lt}(f) = c_{1}x^{a_{1}}y^{b_{1}}$, the \emph{leading term} of $f$.

Since the coefficients are not necessarily in a field, we need to
`re-define' division.

\begin{definition}Let $f$ and $h$ be polynomials in $R$ and $G$ a set of polynomials in $R$, $G = \{g_{1}, g_{2}, ..., g_{n}\}$. We say that $f$
reduces to $h$ modulo the set $G$ in one step, denoted
$$f \overset{G}{\longrightarrow}_{1} h,$$ if and only if $$h = f -
(c_{1}x^{a_{1}}y^{b_{1}}g_{1} + ... +
c_{s}x^{a_{s}}y^{b_{s}}g_{s})$$ for $c_{1}, ..., c_{s} \in (\mathbb{Z}/M\mathbb{Z})$ and
with $\operatorname{lp}(f) = x^{a_{i}}y^{b_{i}}\operatorname{lp}(g_{i})$ for all $i$ such that
$c_{i} \neq 0$ and $\operatorname{lt}(f) = c_{1}x^{a_{1}}y^{b_{1}}\operatorname{lt}(g_{1}) + ... +
c_{s}x^{a_{s}}y^{b_{s}}\operatorname{lt}(g_{s})$.\end{definition}

\begin{definition} Let $f$, $h$ and $f_{1}, f_{2}, ..., f_{s}$ be
polynomials in $R$, with $f_{i} \neq 0 \ \forall \ 1 \leq i \leq s$,
and let $F = \{f_{1}, f_{2}, ..., f_{s}\}$. We say that $f$ reduces
to $h$ modulo $F$, denoted $$\xymatrix{f\ar[r]^{F} & h},$$ if and
only if there exist polynomials $h_{1}, ..., h_{t-1} \in R$ such
that
$$f \overset{F}{\longrightarrow}_{1} h_{1} \overset{F}{\longrightarrow}_{1} h_{2} \overset{F}{\longrightarrow}_{1} ... \overset{F}{\longrightarrow}_{1} h_{t-1} \overset{F}{\longrightarrow}_{1} h.$$ 
We note that if
$$\xymatrix{f\ar[r]^{F} & h},$$ then $f - h \in \langle f_{1}, ...,
f_{s} \rangle$.\end{definition}

We will now give the statement of a theorem (\cite[Theorem
4.1.12]{Adams}) which basically serves as the definition for a
Gr\"obner Basis. We need to state first that \emph{the leading term
ideal} of an ideal $V$ of our ring $R$, denoted by $LT(V)$, is defined
as:
\begin{center}$LT(V) = \langle \{\operatorname{lt}(v) : v \in V\} \rangle.$
\end{center}

\begin{theorem} Let $V$ be an ideal of $R$ and let $G =
\{g_{1}, ..., g_{n}\}$ be a set of non-zero polynomials in $V$. The
following are equivalent:

(i) $LT(G) = LT(V)$.

(ii) For any polynomial $f \in R$ we have \begin{center} $f \in V$
if and only if $$\xymatrix{f\ar[r]^{G} & 0}$$\end{center}

(iii) For all $f \in V$, $f = h_{1}g_{1}+...+h_{n}g_{n}$ for some
polynomials $h_{1}, ..., h_{n} \in R$, such that $\operatorname{lp}(f) = max_{1
\leq i \leq n} (\operatorname{lp}(h_{i})\operatorname{lp}(g_{i}))$. \ $\Box$ \end{theorem}

\begin{definition} A set $G$ of non-zero polynomials contained in an
ideal $V$ of our ring $R$ is called a Gr\"obner basis for $V$ if and
only if $G$ satisfies any one of the three equivalent conditions of
Theorem 4.6 above. Obviously $G$ is a Gr\"obner basis for $\langle G
\rangle$.\end{definition}

The Noetherian property of the ring $R$ and Theorem 4.6 above, yield
the following Theorem (\cite[Corollary 4.1.17]{Adams}):

\begin{theorem} Let $V \subseteq R$ be a non-zero ideal. Then
$V$ has a finite Gr\"obner Basis. \ $\Box$ \end{theorem}

We compute in MATHEMATICA a Gr\"obner basis for our ideal $J$ of the ring $R$, which we denote by $G_{J}$.
We see that the order of
$B[M]_{d_{1},d_{2}}^{\perp}$ is the order of the quotient $ I_{d} /  \langle G_{J} \rangle$.

In the last step of the algorithm we will also need to compute the
annihilator of some ideal $\langle G_{J} \rangle$ over the finite
ring $R_{d_{1},d_{2}}$. For this we follow
the method outlined in \cite[Proposition 4.3.11]{Adams} and we
calculate the ideal quotient
\begin{center}$T : \langle G_{J} \rangle = \{ f \in R : f \langle
G_{J} \rangle \subseteq T \}$\end{center} where $T = \langle
(x^{l^{a_{1}}})^{d_{1}}-1,(y^{l^{a_{2}}})^{d_{2}}-1 \rangle$ is the zero ideal of $R_{d_{1},d_{2}}$.
We see that
\begin{center}$\operatorname{Ann}_{R_{d_{1},d_{2}}}(\langle G_{J} \rangle) = T :
\langle G_{J} \rangle$.\end{center}

\section{The Algorithm} 

\label{Algorithm}

\subsection{Step 1} \label{Step1}

Fix distinct odd primes $p$ and $q$ and an odd prime $l$. The
product $pq$ is the conductor of the field $\mathbb{Q}(\zeta_{pq})^{+}$ whose
class number $h^{+}$ we want to calculate and $M = l$ is the prime
that we check to see if it divides $h^{+}$. Factor $x^{m_{1}}-1$ and
$y^{m_{2}}-1$ into irreducibles in $\mathbb{Z}/l\mathbb{Z}$ where, as above,
$\operatorname{gcd}(m_{i},l)$=1 for $i$=1,2 and $m_{1}l^{a_{1}}=p-1$ and
$m_{2}l^{a_{2}}=q-1$. As before, let $(f_{1},f_{2})$ be a pair of
degrees of irreducible polynomials $\phi_{x}$, $\phi_{y}$
respectively, which appear in the factorization of $\mathbb{Z}[\tilde{G}]$.
Let $d_{1} = \operatorname{gcd}(p-1, l^{f1} - 1)$ and $d_{2} = \operatorname{gcd}(q-1, l^{f_{2}}
-1)$. For various primes $r$ with $r\equiv \pm1~({\rm mod}~pq)$ and $r\equiv 1~({\rm mod}~2l)$ we calculate the Frobenius elements $f_{\Re}$
as in $(iii)$, Section~\ref{Decomposition}. Let $J_{0}$ denote the zero ideal of $R_{d_{1},d_{2}}$
together with the conjugation relation $cnj$. We pick several
Frobenius polynomials ${f_{\Re}}_{i}$ that we calculated above and
we let $J_{i} = J_{i-1} + ({f_{\Re}}_{i})$. This ascending chain of
ideals will computationally stabilize at some ideal $J^{l}$
$\subseteq$ $I_{d}$. If $J^{l}$ happens to
equal the whole ideal $I_{d}$, then the module $B[l]_{d_{1},d_{2}}$ is trivial. If however, for
some pair of degrees $(f_{1},f_{2})$ we have a strict inclusion
$J^{l}$ $\subset$ $I_{d}$ then the corresponding
$B[l]_{d_{1},d_{2}}^{\perp}$ is not trivial, if $J^{l}$ has indeed
stabilized at the correct ideal $J$. Hence we believe that $l$
divides the index [$E$:$H$].

As expected, in most cases the ideal $J^{l}$ is the whole ideal $I_{d}$ and so we do not continue to steps 2 and 3 for
this prime $l$. When we do get a non-trivial quotient
$I_{d}/J^{l}$ for some $l$ however, we do not proceed to the
next step right away but we follow first the procedure outlined
right after the proof of Proposition~\ref{prop1}. 

\subsection{Step 2} \label{Step2}

In this step we repeat the procedure of step 1 but with higher
powers of $l$, i.e. for M =  $l^{2}$, $l^{3}$, etc, and only for
those primes which `passed' step 1. The coefficients of the
Frobenius polynomials $f_{\Re}$ now lie in $(\mathbb{Z}/M\mathbb{Z})$ and we have to
make sure that the primes $r$ satisfy $r\equiv 1~({\rm mod}~2M)$ for
the specific power $M$ of $l$. As before, $R_{d_{1},d_{2}}$ =
$(\mathbb{Z}/M\mathbb{Z})[x,y]/((x^{l^{a_{1}}})^{d_{1}}-1,(y^{l^{a_{2}}})^{d_{2}}-1)$,
and $I_{d}$ is the ideal $(x-1,y-1)$ in $R_{d_{1},d_{2}}$. As in Step 1, for each
$M$ we have that the sequence of ideals $$J_{0} \subset J_{1}
\subset ... \subset J_{i} \subset ...$$ will stabilize at some ideal
$J^{M}$ and the order of the quotients $$..., I_{d}/J^{M}, I_{d}/J^{lM}, ...$$ is non-decreasing. Since
$B[M]_{d_{1},d_{2}}^{\perp}$ is finite and its order is bounded
above by $|B_{d_{1},d_{2}}|$ which is finite and independent of $M$,
the orders of the quotients $I_{d}/J^{M}$ will have to stabilize. We
will have that for some power $M$ of $l$, $|I_{d}/J^{lM}|$ =
$|I_{d}/J^{M}|$ hence $I_{d}/J^{lM} \cong I_{d}/J^{M}$. Therefore
$M$ annihilates $I_{d}/J^{lM}$ as well as its quotient 
$\Big{(} I_{d}/J^{lM} \Big{)} / \langle f_{\Re}(x,y) : \Re \in S \rangle$ $\cong B[lM]_{d_{1},d_{2}}^{\perp}$. This implies
that $M(B[lM]_{d_{1},d_{2}}^{dual}) = 0$ which gives
$M(B[lM]_{d_{1},d_{2}}) = 0$ since as finite abelian groups,
$B[lM]_{d_{1},d_{2}}^{dual}$ and $B[lM]_{d_{1},d_{2}}$ are
isomorphic. Therefore $B[lM]_{d_{1},d_{2}} = B[M]_{d_{1},d_{2}}$ and
$$|MB_{d_{1},d_{2}}| = |B_{d_{1},d_{2}}/B[M]_{d_{1},d_{2}}| =
|B_{d_{1},d_{2}}|/|B[Ml]_{d_{1},d_{2}}| = |lMB_{d_{1},d_{2}}|.$$
Therefore $(MB_{d_{1},d_{2}})/l(MB_{d_{1},d_{2}}) = 0$ and by
Nakayama's Lemma, $MB_{d_{1},d_{2}} = 0$. Again, since
$B_{d_{1},d_{2}}$ and $B_{d_{1},d_{2}}^{dual}$ are isomorphic as
finite abelian groups, we obtain $MB_{d_{1},d_{2}}^{\perp} = 0$. Hence the map $g: I_{d} /J^{M} \rightarrow B_{d_{1},d_{2}}^{\perp}$ is surjective.

Let $GCD(P_{g,h})$ denote the greatest common divisor of all the $P_{g,h}$ that were defined in the last paragraph of Section~\ref{eta}. We need to state here that if $|I_{d}/J^{M}|$ = $|GCD(P_{g,h})|_{l}$,  then we do not have to procced to Step 3 for this prime since 
$|I_{d}/J^{M}| \geq |B_{d_{1},d_{2}}|_{l} \geq |GCD(P_{g,h})|_{l}$ and no power of $l$ divides $h^{+}$.

\subsection{Step 3} \label{Step3}

In the third and last step we determine the structure and hence the
order of the module $B_{d_{1},d_{2}}^{\perp}$, by showing that the
surjective map $g: I_{d} /J^{M} \rightarrow B_{d_{1},d_{2}}^{\perp}$ is actually an isomorphism.

Let $M$ be as in step 2, i.e. the power of $l$ which annihilates
$B_{d_{1},d_{2}}^{\perp}$. Consider the exact sequence, where $\psi'$ raises an element to its $M$-th power:
$$ \xymatrix{ 0 \ar[r] & B[M]\ar[r]^{\psi' \phantom{xyz}} & H/\pm
H^{M} \ar[r] & H /\pm E^{M} \ar[r] & 0 }$$ Since $M$ annihilates $B_{d_{1},d_{2}}^{\perp} \cong
Hom_{R_{d_{1},d_{2}}}(B_{d_{1},d_{2}}, R_{d_{1},d_{2}})$ that
implies that $M$ also annihilates $B_{d_{1},d_{2}}$. Therefore, we obtain the following exact
sequence of $R_{d_{1},d_{2}}$-modules $$ \xymatrix{ 0 \ar[r] &
B_{d_{1},d_{2}} \ar[r]^{\psi \phantom{xyz}} & (H /\pm H^{M})_{d_{1},d_{2}}
\ar[r] & (H /\pm E^{M})_{d_{1},d_{2}} \ar[r] & 0 }$$
where the generator $\eta_{d_{1},d_{2}}$ of the unit groups  in the sequence above, is the unit $\eta$ with the norm map
$N_{d} = \frac{(x^{p-1}-1)(y^{q-1}-1)}{(x^{l^{a_{1}d_{1}}}-1)(y^{l^{a_{2}d_{2}}}-1)}$ applied to it.

From the surjection $g : I_{d}/J^{M} \to B_{d_{1},d_{2}}^{\perp}$ that we
established from step 2 we obtain the injection $$\Psi :
B_{d_{1},d_{2}} \hookrightarrow (I_{d}/J^{M})^{\perp}.$$
We have $(I_{d}/J^{M})^{\perp} = Hom_{R_{d_{1},d_{2}}}(I_{d}/J^{M},R_{d_{1},d_{2}}) \cong R_{d_{1},d_{2}}/\operatorname{Ann}_{R_{d_{1},d_{2}}}(I_{d}/J^{M})$, where $\operatorname{Ann}_{R_{d_{1},d_{2}}}(I_{d}/J^{M}) = \operatorname{Ann}_{R_{d_{1},d_{2}}}(I_{d}+J^{M}/J^{M}) = (J^{M} : I_{d})$, which is the ideal quotient already discussed in Section~\ref{Groebner}. Therefore, $(I_{d}/J^{M}) \cong \operatorname{Ann}_{R_{d_{1},d_{2}}}((J^{M} : I_{d}))$. One can think of the ideal quotient $(J^{M} : I_{d})$, which we denote by $\overline{J^{M}}$, as the ideal $J^{M}$  modulo its $I_{d}$-part.
Assume now that $\operatorname{Ann}_{R_{d_{1},d_{2}}}(\overline{J^{M}})$ annihilates
$(H /\pm H^{M})_{d_{1},d_{2}}/\psi(B_{d_{1},d_{2}})$. Then
$\operatorname{Ann}_{R_{d_{1},d_{2}}}(\overline{J^{M}}) \subseteq
\psi(B_{d_{1},d_{2}}).$ But now we have that
$$|\operatorname{Ann}_{R_{d_{1},d_{2}}}(\overline{J^{M}})| \leq |\psi(B_{d_{1},d_{2}})| =
|B_{d_{1},d_{2}}| = |\Psi(B_{d_{1},d_{2}})| \leq
|\operatorname{Ann}_{R_{d_{1},d_{2}}}(\overline{J^{M}})|.$$ Therefore we have that the
orders of $I_{d}/J^{M}$ and $B_{d_{1},d_{2}}^{\perp}$ are equal and
hence $g$ is an isomorphism. Hence, if we show that
$\operatorname{Ann}_{R_{d_{1},d_{2}}}(\overline{J^{M}})$ annihilates $(H/\pm
E^{M})_{d_{1},d_{2}}$, from the second exact sequence above we will have proved that $g$ is an
isomorphism.

To find the annihilator $\operatorname{Ann}_{R_{d_{1},d_{2}}}(\overline{J^{M}})$ we first compute the ideal quotient $\overline{J^{M}} = (J^{M} : I_{d})$ and then the ideal quotient $(T : \overline{J^{M}} )$ of Section~\ref{Groebner}. For each generator of $\operatorname{Ann}_{R_{d_{1},d_{2}}}(\overline{J^{M}})$, we need to apply a lift $h(x,y)$ $\in \mathbb{Z}[x,y]$ of this generator to the unit $\eta_{d_{1},d_{2}}$ $\in$ $(H/\pm E^{M})_{d_{1},d_{2}}$. If $\eta_{d_{1},d_{2}}^{h(x,y)}$ is an $M$-th power
of a unit in $E$ then we are done. To see whether it is an $M$-th
power we follow a method similar to the one in Gras and Gras
\cite{Gras2} that we also mentioned in the Introduction. We
reformulate here the main proposition from \cite{Gras2} in order to
make it applicable to our case and we prove it again, only for the
case that $l$ is odd since we only calculate the odd $l$-parts of
$h^{+}$.

We denote by $\eta_{d}^{h}$ the unit $\eta_{d_{1},d_{2}}^{h(x,y)}$
that we already described above and by $G_{d}$ the quotient of $G$
containing the coset representatives of the embeddings in $G$, which
map $\zeta_{p}$ to $\zeta_{p}^{g^{i}}$ and $\zeta_{q}$ to
$\zeta_{q}^{h^{j}}$, for $1 \leq i \leq l^{a_{1}}d_{1}$ and $1  \leq j
\leq l^{a_{2}}d_{2}$. 

\begin{proposition} \label{prop2} Let $M$ be a fixed power of an odd prime $l$ as
above and consider the polynomial $$P(X) = \prod_{a \in G_{d}}
(X-(a(\eta_{d}^{h}))^{1/M})$$ where $(a(\eta_{d}^{h}))^{1/M}$
denotes the real $M$-th root of $a(\eta_{d}^{h})$. If $P$ has
coefficients in $\mathbb{Z}$ then $\eta_{d}^{h}$ is an $M$-th power in
$\mathbb{Q}(\zeta_{pq})^{+}$.\end{proposition}

\emph{Proof:} Let $N$ be the largest power of $l$ for which the unit
$(\eta_{d}^{h})^{1/N}$ lies in $\mathbb{Q}(\zeta_{pq})^{+}$. If $M = N$ then
we are done so we assume $N < M$. Then $(\eta_{d}^{h})^{1/N}$ is not
an element of $(\mathbb{Q}(\zeta_{pq})^{+})^{l}$ and therefore by
\cite[Chapter VIII, Theorem 16]{Lang} we have that the polynomial
\begin{center}$T(X) = X^{M/N} - (\eta_{d}^{h})^{1/N}$\end{center} is
irreducible in $\mathbb{Q}(\zeta_{pq})^{+}$. Since $M/N \geq 3$, $T(X)$ has
at least one complex root. Therefore $(\eta_{d}^{h})^{1/M}$ has at
least one Galois conjugate that is not real. But $P(X) \in \mathbb{Z}[X]$
implies that the Galois conjugates are roots of $P(X)$ which are
real. Therefore we have a contradiction. \ $\Box$

Most of the times the $P(X)$'s are very large polynomials with huge coefficients. We can prove that they are integers by applying a method outlined in Schoof \cite{Schoof1}, which requires that we round off the coefficients of $P(X)$ and then show that this new polynomial divides $P(X^{M})$.

\section{Examples} \label{examples}

\subsection{The field of conductor 469 = 7 $\cdot$ 67}

We confirm Hakkarainen's \cite{Hak} result that 3 is the only odd prime $<$ 10000 which divides $h^{+}$. He
only obtained however a $3^{1}$ dividing $h_{\xi}$, whereas our
results show that that the 3-part of $h^{+}$ has order $3^{2}$.

Let $r$ =3, $p$=7, $q$=67 and $\mathbb{Q}(\zeta_{pq})^{+}$ be the real
cyclotomic field of conductor $pq$=469. We first compute the factor
$P_{(g,h)}$ for all pairs of primitive roots $(g,h)$ and then their
greatest common divisor $GCD(P_{(g,h)})$. From the calculations we have that
$GCD(P_{(g,h)}) = 2^{32}$ and so we see that it is best to run the
test with the pair $(g',h') = (3 (mod 7), 7 (mod 67))$ for which
$P_{(g',h')}$ has the smallest number of factors. In particular,
$P_{(g',h')} = 2^{98}\cdot17^{2}$. Next, we decompose the group ring
$\mathbb{Z}[G]$ as we show in Section~\ref{Decomposition}. We have that $x^{p-1}-1$ =
$(x^{3})^{2}-1$ and $y^{q-1}-1$ = $(y^{3})^{22}-1$ and factoring into irreducibles in $\mathbb{Z}/3\mathbb{Z}$ gives $x^{2}-1 = (x+1) (x+2)$ and $y^{22}-1 = (y+1) (y+2) (y^5+2y^{3}+y^{2}+2y+2) (y^{5}+2y^{3}+2y^{2}+2y+1)(y^{5}+2y^{4}+2y^{3}+2y^{2}+1) (y^{5}+y^{4}+2y^{3}+y^{2}+2)$
and so we run step 1 for all possible degrees $d_{1}$ and $d_{2}$ which in this case are $d_{1} = 2$ and $d_{2} = 2$ and $22$.
Step 1 gave the primes 2, 3 and 17 to be the only primes $<$ 10000
that are possible divisors of the index. Since we chose not to
calculate the 2-part of $h^{+}$, the only primes we have to consider
are 3 and 17. Before proceeding  to step 2 however, we run
step 1 again for the prime 17 because it did appear as a factor of
$P_{(g',h')}$ but not of $GCD(P_{(g,h)})$ and therefore it is
possible that it might only divide $P_{(g',h')}$ and not $h^{+}$.
The pair $(g_{0},h_{0}) = (5 (mod7), 7 (mod67))$ does not have 17 as
a factor of $P_{(g_{0},h_{0})}$ and step 1 for 17 with this pair of
primitive roots only gives trivial Jordan-H\"older factors.
Therefore we proceed to the next steps only for the prime 3. In step 2 we repeat the same procedure as in step 1 but with higher
powers of 3. Below we show the Frobenius polynomials obtained for $M =3, 3^{2}$ and
$3^{3}$ for the pair of degrees $(d_{1},d_{2})$ = (2,2), the ideals
$J^{M}$ at which the ideals $J_{i}$ stabilize and the order of the
quotients $|I_{d}/J^{M}|$. 

For $M = 3$, $J^{M} = (y^{2} - 1, y - x) \equiv ((y+1)(y-1), (y-1) - (x-1))$ with coefficients in
$\mathbb{Z}/3\mathbb{Z}$. The relation $(y-1) - (x-1)$ implies that $(y-1)$ and $(x-1)$ are equivalent. From the second polynomial in $J^{M}$ we see that the two generators
of the augmentation ideal $I_{d}$ become equivalent in
$I_{d}/J^{M}$. From the first one we have that $y(y-1) \equiv
-(y-1)$ in $J^{M}$ therefore we can only have constants in front of
the only generator of $I_{d}/J^{M}$. Since we are in $\mathbb{Z}/3\mathbb{Z}$ we have that $|I_{d}/J^{M}|$ = 3.

For $M = 9$, $J^{M} = (y^{2}-3y+2, 3-x-2y) = ((y-1)(y-2), -2(y-1)-(x-1))$ with coefficients in
$\mathbb{Z}/9\mathbb{Z}$. The same reasoning as above for the ideal $J^{3}$ applies here as
well and we have that $|I_{d}/J^{M}| = 3^{2}$. Since $|I_{d}/J^{3}|$
is strictly smaller than  $|I_{d}/J^{3^{2}}|$ we
need to continue as above with $M = 3^{3}$.

For $M= 27$, $J^{M} = (9(y-1), 2-3y+y^{2}, 3-x-2y)$ with coefficients in $\mathbb{Z}/27\mathbb{Z}$.
We see here that $J^{3^{3}}$ is generated by the same polynomials as
$J^{3^{2}}$ but it has the extra polynomial $9(y-1)$ which reduces
the number of constants to 9 instead of 27. Therefore
$|I_{d}/J^{3^{2}}|$ = $|I_{d}/J^{3^{3}}|$ = 9 and so, as
expected, the orders of these quotients stabilize with $M = 3^{2}$. The factors $\phi_{x}$, $\phi_{y}$ contained in the ideals $J$ above are $\phi_{x}(x) = x+1$ and $\phi_{y}(y)=y-2$. These are the only two factors that gave a non-trivial quotient $I_{d}/J^{M}$.

For the pair of degrees $(d_{1},d_{2})$ = (2,22) the Frobenius
polynomials give exactly the same ideals $J^{M}$ as above and
therefore we have the same two factors $\phi_{x}$ and $\phi_{y}$. Hence, we only need to consider the case for $(d_{1},d_{2})$ =
(2,2). 

\begin{center}
\tiny
\begin{tabular}{|c| p{11cm} |}
\hline
$l$ & The Frobenius Maps for $M$ = 3 \\
\hline \hline
$l_{1}=7521823$ & $f_{\Re_{1}} = (y^{5}+y^{4}+2y^{3}+2y^{2}+2y+2)x^{5}+(2y^{5}+2y^{4}+2y^{3}+2y^{2}+y+2)x^{4}+(2y^{5}+y^{4}+y^{2}+1)x^{3}+(2y^{5}+2y^{4}+2y^{3}+2y^{2}+y+2)x^{2}+(y^{5}+y^{4}+y^{3}+2y+1)x+(y^{5}+y^{4}+y^{3}+2y)$ \\
\hline
$l_{2} = 8889427$ & $f_{\Re_{2}} = (2y^{5}+2y^{4}+y^{3}+2y^{2}+y+2)x^{5}+(2y^{5}+2y^{3}+2y^{2}+2y)x^{4}+(2y^{5}+2y^{4}+y^{3}+ 2y+2)x^{3}+(y^{5}+2y^{3}+2y^{2}+2y)x^{2}+(2y^{4}+y)x+(y^{5}+2y^{4}+2y^{3}+2y^{2}+y)$ \\
\hline
$l_{3} = 9573229$ & $f_{\Re_{3}} = (y^4 + 2y^3 + 2y + 1)x^5 + (y^4 + y^3 + y^2)x^4 + (y^5 + 2y^2 + y + 2)x^3 + (2y^4 + 2y^3)x^2 + (2y^5 + 2y^4 + y^3 + y + 2)x + (y^5 + y^3 + 1)$ \\
\hline
$l_{4} = 10257031$ & $f_{\Re_{4}} = (y^5 + y + 2)x^5 + (2y^5 + 2y^4 + 2y^3 + 2y^2)x^4 + (2y^4 + 2y^3 + y^2 + 2y)x^3 + (y^5 + y^2 + y + 1)x^2 + (2y^5 + 2y^4 + y^2 + 2y + 2)x + (2y^5 + y^4 + 2y^3 + y^2 + y)$ \\
\hline
$l_{5} = 20514061$ & $f_{\Re_{5}} = (2y^5 + y^3 + y^2 + 2y + 1)x^5 + (2y^4 + y^3 + y^2 + 2y + 2)x^4 + (2y^5 + y^4 + 2y^3 + y + 2)x^3 + (y^2 + y)x^2 + (2y^5 + y^2 + 2y + 1)x + (y^5 + 2y^3 + y^2 + y)$ \\
\hline
$l_{6} = 22565467$ &  $f_{\Re_{6}} = (2y^{4} + y^{3} + y^{2} + 2)x^{5} + (2y^{5} + 2y^{3} + y + 1)x^{4} + (y^{5} + y^{4} + y^{3} + y^{2} + 2y + 1)x^{3} + (2y^{4} + 2y^{3} + y^{2} + y + 2)x^{2} + (y^{5} + 2y^{4}+ 2y^{3} + y^{2} + 2)x + (2y^{5} + 2y^{4} + y^{3} + y^{2} + 1)$ \\
\hline \hline
$l$ & The Frobenius Maps for $M = 3^{2}$\\
\hline \hline
$l_{1}=7521823$ & $f_{\Re_{1}} = (4y^{5} + 7y^{4} + 5y^{3} + 8y^{2} + 5y + 2)x^{5} + (5y^{5} + 8y^{4} + 5y^{3} + 2y^{2} + 4y + 8)x^{4} + (8y^{5} + 4y^{4} + 6y^{3} + 7y^{2} + 6y + 4)x^{3} + (2y^{5} + 8y^{4} + 2y^{3} +
2y^{2} + y + 8)x^{2} + (y^{5} + 7y^{4} + y^{3} + 3y^{2} + 8y + 7)x + (y^{5} + 4y^{4} + 7y^{3} + 3y^{2} + 2y + 6)$ \\
\hline
$l_{2} = 8889427$ & $f_{\Re_{2}} = (2y^{5} + 2y^{4} + y^{3} + 2y^{2} + 4y + 5)x^{5} + (2y^{5} + 3y^{4} + 8y^{3} + 2y^{2} + 5y + 3)x^{4} + (2y^{5} + 8y^{4} + y^{3} + 5y + 5)x^{3} + (7y^{5} + 6y^{4} + 5y^{3} + 2y^{2} + 2y +
6)x^{2} + (6y^{5} + 2y^{4} + 3y^{2} + y)x + (y^{5} + 2y^{4} + 5y^{3} + 2y^{2} + y + 6)$\\
\hline
\end{tabular}
\end{center}

\begin{center}
\tiny
\begin{tabular}{|c| p{11cm} |}
\hline
$l_{3} = 9573229$ & $f_{\Re_{3}} = (4y^{4} + 8y^{3} + 5y + 1)x^{5} + (3y^{5} + 7y^{4} + 7y^{3} + 7y^{2} + 6y)x^{4} + (y^{5} + 6y^{3} + 5y^{2} + y + 2)x^{3} +(6y^{5} + 5y^{4} + 5y^{3} + 3y^{2} + 3y + 3)x^{2} +
(5y^{5} + 5y^{4} + y^{3} + 7y + 8)x + (4y^{5}+ 3y^{4} + 7y^{3} + 3y^{2} + 6y + 7)$ \\
\hline
$l_{4} = 10257031$ & $f_{\Re_{4}} = (y^{5} + 3y^{4} + y + 2)x^{5} + (2y^{5} + 5y^{4} + 2y^{3} + 8y^{2})x^{4} + (6y^{5} + 2y^{4} + 5y^{3} + y^{2} + 5y)x^{3} +(4y^{5} + 7y^{2} + y + 4)x^{2} + (8y^{5} + 5y^{4} + 6y^{3} + y^{2}
+ 2y + 8)x + (8y^{5} + 7y^{4} + 8y^{3} + 4y^{2}+ 4y + 6)$ \\
\hline
$l_{5} = 20514061$ & $f_{\Re_{5}} =  (5y^{5} + 3y^{4} + 7y^{3} + 7y^{2} + 2y + 4)x^{5} + (5y^{4} + 7y^{3} + 4y^{2} + 2y + 2)x^{4} + (5y^{5} + 4y^{4} + 8y^{3} + 3y^{2} + y + 2)x^{3} + (3y^{5} + 6y^{4} + 3y^{3} + y^{2} + y +
3)x^{2} + (8y^{5} + 6y^{4} + 3y^{3} + 7y^{2} + 8y + 7)x + (4y^{5} + 3y^{4} + 8y^{3} + 7y^{2} + 7y + 6)$ \\
\hline
$l_{6} = 22565467$ & $f_{\Re_{6}} = (3y^{5} + 5y^{4} + y^{3} + y^{2} + 3y + 2)x^{5} + (5y^{5} + 6y^{4} + 8y^{3} + 7y + 7)x^{4} + (y^{5} + 7y^{4} + 7y^{3} + 4y^{2} + 8y + 1)x^{3} + (3y^{5} + 5y^{4} + 2y^{3} + 7y^{2} + 7y +
8)x^{2} + (7y^{5} + 2y^{4} + 8y^{3} + 7y^{2} + 5)x + (2y^{5} + 5y^{4} + y^{3} + y^{2} + 7)$\\
\hline \hline
$l$ & The Frobenius Maps for $M = 3^{3}$ \\
\hline \hline
$l_{1} = 7521823$ & $f_{\Re_{1}} = (13y^{5} + 7y^{4} + 14y^{3} + 26y^{2} + 14y + 20)x^{5} + (5y^{5} + 17y^{4} + 5y^{3} + 11y^{2} + 13y + 26)x^{4} + (17y^{5} + 4y^{4} + 24y^{3} + 16y^{2} + 6y + 4)x^{3} + (20y^{5} +
8y^{4} + 11y^{3} + 11y^{2} + y + 26)x^{2} + (y^{5} + 25y^{4} + 19y^{3} + 21y^{2} + 26y + 7)x + (10y^{5} + 22y^{4} + 25y^{3} + 3y^{2} + 2y + 6)$\\
\hline
$l_{2} = 8889427$ & $f_{\Re_{2}} = (20y^{5} + 2y^{4} + y^{3} + 11y^{2} + 22y + 23)x^{5} + (2y^{5} + 12y^{4} + 26y^{3} + 11y^{2} + 14y + 21)x^{4} + (2y^{5} + 26y^{4} + y^{3} + 18y^{2} + 14y + 23)x^{3} + (16y^{5} + 15y^{4} +
14y^{3} + 20y^{2} + 11y + 15)x^{2} + (6y^{5} + 11y^{4} + 21y^{2} + 19y)x + (10y^{5} + 11y^{4} + 23y^{3} + 11y^{2} + 10y + 24)$\\
\hline
$l_{3} = 9573229$ & $f_{\Re_{3}} = (9y^{5} + 4y^{4} + 17y^{3} + 18y^{2} + 5y + 19)x^{5} + (12y^{5} + 16y^{4} + 7y^{3} + 7y^{2} + 15y)x^{4} + (10y^{5} + 24y^{3} + 14y^{2} + 10y + 20)x^{3} + (24y^{5} + 5y^{4} + 5y^{3} +
3y^{2} + 12y + 12)x^{2} + (5y^{5} + 5y^{4} + 19y^{3} + 9y^{2} + 7y + 26)x + (13y^{5} + 21y^{4} + 25y^{3} + 12y^{2} + 6y + 16)$ \\
\hline
$l_{4} = 10257031$ & $f_{\Re_{4}} = (10y^{5} + 3y^{4} + 18y^{3} + 19y + 2)x^{5} + (20y^{5} + 14y^{4} + 2y^{3} + 17y^{2} + 9y + 9)x^{4} + (24y^{5} + 2y^{4} + 23y^{3} + 10y^{2} + 14y + 9)x^{3} + (22y^{5} + 9y^{4} +
9y^{3} + 7y^{2} + 10y + 22)x^{2} + (26y^{5} + 5y^{4} + 15y^{3} + y^{2} + 2y + 26)x + (17y^{5} + 16y^{4} + 26y^{3} + 4y^{2} + 22y + 15)$ \\
\hline
$l_{5} = 20514061$ & $f_{\Re_{5}} =  (14y^{5} + 3y^{4} + 25y^{3} + 16y^{2} + 11y + 22)x^{5} + (18y^{5} + 23y^{4} + 16y^{3} + 22y^{2} + 2y + 20)x^{4} + (23y^{5} + 22y^{4} + 17y^{3} + 21y^{2} + 10y + 20)x^{3} + (21y^{5} +
6y^{4} + 12y^{3} + 19y^{2} + y + 3)x^{2} + (17y^{5} + 15y^{4} + 21y^{3} + 7y^{2} + 26y + 16)x + (13y^{5} + 3y^{4} + 26y^{3} + 25y^{2} + 7y + 24)$ \\
\hline
$l_{6} = 22565467$ & $f_{\Re_{6}} = (12y^{5} + 23y^{4} + 19y^{3} + 10y^{2} + 12y + 11)x^{5} + (14y^{5} + 24y^{4} + 8y^{3} + 25y + 7)x^{4} + (10y^{5} + 16y^{4} + 7y^{3} + 4y^{2} + 8y + 1)x^{3} + (21y^{5} + 5y^{4} +
20y^{3} + 16y^{2} + 16y + 8)x^{2} + (25y^{5} + 2y^{4} + 8y^{3} + 16y^{2} + 5)x + (20y^{5} + 14y^{4} + y^{3} + y^{2} + 16)$\\
\hline \hline
\end{tabular}
\end{center}

\normalsize 

We now proceed to step 3 of the algorithm where we prove that
$I_{d}/J^{M}$ is isomorphic to $B_{d_{1},d_{2}}^{\perp}$. The computations for the ideal $\overline{J^{M}}$ gave us the basis $(-2 + y, 1 + x)$ in $\mathbb{Z}/9\mathbb{Z}$. Since the degree of x and y in $R_{d_{1},d_{2}}$ is $3\cdot2=6$, we compute the annihilator of $(-2 + y^{3}, 1 + x^{3})$ in $(\mathbb{Z}/9\mathbb{Z})[x,y]/(x^{6}-1,y^{6}-1)$. We found the
following polynomial to be the generator of
$\operatorname{Ann}_{R_{d_{1},d_{2}}}(\overline{J^{M}})$: $h(x,y) = 3 - 3 x^3 - 3 y^3 + 3 x^3 y^3$.
Factoring $h(x,y)$ in $\mathbb{Z}[x,y]$ we find $h(x,y)=3 (-1 + x) (1 + x + x^2) (-1 + y) (1 + y + y^2) = 3(x-1)(y-1)\Phi_{3}(x)\Phi_{3}(y)$ where $\Phi_{k}$ is the $k$-th cyclotomic polynomial. Therefore, we apply to $\eta$ the norm map  $\frac{(x^{6}-1)(y^{66}-1)}{(x^{6/3}-1)(y^{6/3}-1)}$ instead of the map  $\frac{(x^{6}-1)(y^{66}-1)}{(x^{6}-1)(y^{6}-1)}$ and then the annihilator $h'(x,y)=3 (-1 + x) (-1 + y)$. The polynomials $P(x)$ and $P(x^{M})$ of Proposition~\ref{prop2} were calculated with a precision of 500 and are shown in the table below. Finally, we showed that $P(x)$ divides $P(x^{M})$ hence proving rigorously that $3^{2} || h^{+}$. 

\begin{center}
\tiny
\begin{tabular}{|c| p{11.2cm} |}
\hline \hline
$P(x)$ & $x^{4} - 35667454 \cdot x^{3} + 318041818710531 \cdot x^{2} - 35667454 \cdot x + 1$\\
\hline
$P(x^{M})$ & $x^{36} - 364929542762806942594907901654249278525439344697663012299174707204 \cdot x^{27} + 33293392795267835243258623959180895487795677296162956508170492359406402$ \\
 & $192775112912608077373493763985920516781456745581649782374406 \cdot x^{18} -$ \\
 & $364929542762806942594907901654249278525439344697663012299174707204 \cdot x^{9} + 1$\\
\hline \hline
\end{tabular}
\end{center}

\subsection{Step 3 for the field of conductor 1477 = 7 $\cdot$ 211}

We found that the only primes $l < 10000$ dividing $h^{+}$ are $7$ and $11$. In this example we will show our work for $l = 7$ where there are more than one pairs of $(\phi_{x},\phi_{y})$ contained in $J^{M}$.  

Step 2 of our Algorithm showed that the orders of the quotients $I_{d}/J^{M}$ stabilize at $M = 7$ with $J^{M} = (3 + 3 y + y^2, 3 + 6 x + 4 y + x y, 5 + x^3 + y)$ in $\mathbb{Z}/7\mathbb{Z}$ for both possible combinations of $(d_{1},d_{2})$ and we therefore choose the smaller pair $(6,6)$. The pairs of $(\phi_{x},\phi_{y})$ contained in $J^{M}$ are $(x+3,y-1)$, $(x+4,y+4)$, $(x+5,y-1)$. As expected, $|B[M]_{d_{1},d_{2}}^{\perp}|$ = $\prod_{\phi_{x},\phi_{y}} |B[M]_{\phi_{x},\phi_{y}}^{\perp}|$ = $7^{3}$ (we found $|B[M]_{\phi_{x},\phi_{y}}^{\perp}| = 7$ for each pair $(\phi_{x},\phi_{y})$). We have that $|GCD(P_{g,h})|_{7} = 7^{2}$ hence, after Step 3, we will have proved that $7 || h^{+}$.

The computations for the annihilator gave us only one generator 
\begin{center} \tiny{$h(x,y) =  \Phi_{3}(x) \Phi_{3}(y) \Phi_{21}(y) (4 + 2 x - 4 x^{2} + 5 x^{3} + 2 y^{7} - 6 x y^{7} + 
   5 x^{2} y^{7} - x^{3} y^{7} - 4 y^{14} + 5 x y^{14} - 3 x^{2} y^{14} + 
   2 x^{3} y^{14} + 5 y^{21} - x y^{21} + 2 x^{2} y^{21} + x^{3} y^{21}) = \Phi_{3}(x) \Phi_{3}(y) \Phi_{21}(y) h'(x,y)$} \end{center} 
 
\normalsize
 
From basic properties of cyclotomic polynomials, we can establish the relation \scriptsize{$$\Gamma'_{\delta} = \prod_{k \in D(|\Gamma|) \setminus D(\delta)} \Phi_{k}(\sigma) \phantom {xyz} (*)$$} \normalsize where $\Gamma$ is any finite cyclic group of order $|\Gamma|$, $D(|\Gamma|) = \{k \in N : k | |\Gamma| \}$ and $\Gamma'_{\delta}$ is the sum of all the elements of its subgroup $\Gamma_{\delta} = \langle \sigma^{\delta} \rangle$ of index $\delta$ in $\Gamma$.

Let's recall for a moment that our Galois group $\tilde{G}$ of the extension $\mathbb{Q}(\zeta_{pq})/\mathbb{Q}$ is such that $\tilde{G} \cong \tilde{G_{1}} \times \tilde{G_{2}}$ where $\tilde{G_{i}}$ finite cyclic of order $p-1=6$ and $q-1=270$ in this example, for $i = 1$ and $2$ respectively.  As $(d_{1},d_{2}) = (6,6)$ and $l^{a_2} = 7$ we have the subgroup $\tilde{G_{1}} \times \tilde{G_{2}}$ of $\tilde{G}$ of order $6 \cdot 42$.  As we work in the real subfield $\mathbb{Q}(\zeta_{pq})^{+}$, we can assume that the subgroup $G_{1} \times G_{2}$ of $G = Gal(\mathbb{Q}(\zeta_{pq})^{+}/\mathbb{Q})$ is of order $6 \cdot 21$. We then see that $\Phi_{3}(y) \cdot \Phi_{21}(y) = (G_{2})'_{7}$ which, according to the formula $(*)$ above, is the sum of the elements of the group $(G_{2})_{7}$ of index $7$ in $G_{2}$. The norm map that we therefore need to apply to $\eta$ is $\frac{(x^{6}-1)(y^{210}-1)}{(x^{6/3}-1)(y^{42/7}-1)}$ and then the annihilator $h'(x,y)$. The polynomials $P(x)$ and $P(x^{M})$ of Proposition~\ref{prop2} were calculated with a precision of 1000 and we showed that $P(x)$ divides $P(x^{M})$ hence proving rigorously that $7 || h^{+}$. The polynomial $P(x)$ is shown in the table below whereas $P(x^{M})$ is omitted as it is too long.

\begin{center}
\tiny
\begin{tabular}{|c| p{12cm} |}
\hline \hline
$P(x)$ & $x^{12} - 253285672818085597920117540833320566764\cdot x^{11} + 16038408013727576378675398$\\
 &$205615384849932252547671390045959497056856423999746\cdot x^{10}$\\
 &$ - 7447696110433675817548561818649227038803699459085663820108397789285266813843540$\\
 &$443700\cdot x^{9} + 86911768356572921123499159325706075635679782571876578383262$\\
 &$7442647577671131436873144346703615\cdot x^{8} - 1045502371457459906661385781160012228$\\  
 & $359906832607656705019244404117434162718039857442427623338144074\cdot x^{7} + $\\
 & $31606394380090436053358292646577202064791027624187676425$\\
 &$1660067572877683081961563342938807318368386032296\cdot x^{6} - 332994422221005688150$\\
 &$10667879016823147457384730252190107458282081938117105395253036177107901993728$\\
  &$31287574\cdot x^{5} + 877094806999502083991271352174430122151878622656805163463$\\
 &$0457643479432149168032626544704205804528987979615\cdot x^{4} -$ \\
 &$1944459765899336452214557670670811109932072061509042393676$\\
  &$0599147020507012734321450\cdot x^{3} + 10777026227137095866981035797948$\\
 &$453135069447299390696542871\cdot x^{2} - 207623109700797451167702365014\cdot x + 1$\\
\hline \hline
\end{tabular}
\end{center}

\subsection{Step 3 for the field of conductor 1355 = 5 $\cdot$ 271}

We found that the only prime $l < 10000$ dividing $h^{+}$ is $37$. Step 2 of our Algorithm showed that the orders of the quotients $I_{d}/J^{M}$ stabilize at $M = 37$ with $J^{M} = (9 + 27 y + y^2, 8 + x + 28 y)$ in $\mathbb{Z}/37\mathbb{Z}$ and with $(d_{1},d_{2}) = (4,18)$. The only two factors that gave a non-trivial quotient $I_{d}/J^{M}$ where $\phi_{x}(x) = x+1$ and 
$\phi_{y}(y)=y+28$. The computations for the annihilator gave us only one generator 
 \begin{center} \tiny{$= \Phi_{4}(x) \Phi_{2}(y) \Phi_{6}(y) \Phi_{18}(y) (4 + 33 x + 21 y + 
   16 x y + 27 y^{2} + 10 x y^{2} + 3 y^{3} + 34 x y^{3} + 25 y^{4} + 
   12 x y^{4} + 11 y^{5} + 26 x y^{5} + 30 y^{6} + 7 x y^{6} + 28 y^{7} + 
   9 x y^{7} + 36 y^{8} + x y^{8})$} \end{center} 
 \begin{center} \tiny{$= \Phi_{4}(x) \Phi_{2}(y) \Phi_{6}(y) \Phi_{18}(y) h'(x,y)$} \end{center} 
 
 \normalsize
 
According the formula $(*)$ above, $\Phi_{4}(x) = (G_{1})'_{2}$ and $\Phi_{2}(y) \Phi_{6}(y) \Phi_{18}(y) = (G_{2})'_{9}$ hence, the norm map that we need to apply to $\eta$ is 
$\frac{(x^{4}-1)(y^{270}-1)}{(x^{4/2}-1)(y^{18/9}-1)}$ and then the annihilator $h'(x,y)$. The polynomials $P(x)$ and $P(x^{M})$ of Proposition~\ref{prop2} were calculated with a precision of 7000 and we showed that $P(x)$ divides $P(x^{M})$ hence proving rigorously that $37 || h^{+}$. The polynomial $P(x)$ is shown in the table below whereas $P(x^{M})$ is omitted as it is too long.

\begin{center}
\tiny
\begin{tabular}{|c| p{12cm} |}
\hline \hline
$P(x)$ & $x^4 - 534186444472275956720533076216968091508192072459731400996*x^3 + $\\
 &$71338789364482991009380877708435286461900572928062358768758453333592004560$\\
 &$744265699158031988648292121436237448006*x^2 - 534186444472275956720533076216$\\
 &$968091508192072459731400996*x + 1$\\
 \hline \hline
\end{tabular}
\end{center}
 
\section{Table and Discussion of the Results}

We applied Steps 1 and 2 of the algorithm to real cyclotomic fields of conductor $pq$ $<$ 2000. We tested the divisibility of $|B| = [E : H] = P \cdot h^{+}$ by all primes $l$ $<$ 10000. All the primes appearing in the greatest common divisor of the $P_{(g,h)}$ for all pairs of primitive roots $(g,h)$, came up as possible divisors, as expected. These primes are listed in the column `$GCD$' of Table 1 below. Since we do not calculate the 2-part of $h^{+}$ or the $l$-part for $l$ $>$ 10000, we leave out the powers of 2 as well as the primes $>$ 10000 from the `$GCD$'. Therefore, if a `1' appears in the column `$GCD$' for some field, this means that no odd primes $<$ 10000 divide the greatest common divisor of the various $P_{(g,h)}$. However, there are always powers of 2 in the `$GCD$', as we see from our calculations of the index $[E:H]$ in Section 2. In the column `$l$' we present all other primes that Step 1 gave to be possible divisors of $h^{+}$, besides the ones that already appear in the column `$GCD$'. Step 2 verified for all fields that indeed, for some $M$ we have that $|B[M]_{d_{1},d_{2}}|$ = $|B[lM]_{d_{1},d_{2}}|$. As mentioned right before Section~\ref{Step3} above, we proceed to Step 3 only for those primes $l$ with $|I_{d}/J^{M}| > |GCD(P_{g,h})|_{l}$. Finally, in the column `$Degree$' we list the smallest degrees $(d_{1}l^{a_{1}},d_{2}l^{a_{2}})$ for which the module $B^{\perp}_{d_{1},d_{2}}$ turned out to be non trivial and the column `$\tilde{h}^{+}$' shows the $l$-part of $h^{+}$ for all odd primes $l <$ 10000.

For the primes $l$ not dividing the degree of the extension, our results agree with Hakkarainen \cite{Hak}. For those primes $l$ that divide the degree of the extension, our results completed his results in the sense that we either verified that no higher powers of $l$ divide $h^{+}$, or that there are indeed such higher powers of $l$. We found five such fields where higher powers of $l$ divide $h^{+}$, and we mark them with an asterisk under the column `$\tilde{h}^{+}$' in Table 1 below.

\begin{center}
\tiny
\begin{center}\textbf{Table 1}\end{center}
\begin{tabular}{|p{0.7cm}|p{1.2cm}|p{0.4cm}|p{1.2cm}|p{1.3cm}||p{0.7cm}|p{1.2cm}|p{0.4cm}|p{1.2cm}|p{1.2cm}|}
\hline
$f$ & GCD & $l$ & Degree & $\tilde{h}^{+}$ & $f$ & GCD & $l$ & Degree & $\tilde{h}^{+}$\\
\hline \hline
3$\cdot$107 & 1 & 3 & (2,2) & 3 & 
7$\cdot$61 & 1 & 5 & (2,20) & 5 \\
7$\cdot$67 & 1 & 3 & (6,6) & $3^{2}$ $(\ast)$ &
11$\cdot$43 & $3^{4}\cdot5^{2}\cdot7^{4}$ & - & (10,6) & 3 \\
13$\cdot$37 & $7\cdot19$ & - & (6,18) & 19 & 
19$\cdot$29 & $5$ & - &  (2,4) & 5 \\ 
17$\cdot$37 & $3^{4}\cdot19$ & 5 & (4,4) (2,18)& 5$\cdot$19 &
17$\cdot$41 & $3^{3}\cdot7$ & - & (2,2)& 3 \\
19$\cdot$37 & $3^{16}\cdot5$ & 13 & (6,12) &13$\cdot$37 &
3$\cdot$251 & 1 & 11 & (2,10) &11\\
 &  & 37 & (18,36) &  &
  & & & & \\
  7$\cdot$109 & $3^{4}$ & 13 &(6,12) & 13 & 
19$\cdot$41 & $5^{2}$ & 41 & (2,40) &41 \\
 5$\cdot$157 & $3^{2}\cdot79$ & - &(2,6)& 3 & 
 13$\cdot$61 & $3^{20} \cdot5 \cdot 7$ & 37&(12,12) (6,60)& $3\cdot37 \ (\ast)$ \\
 19$\cdot$43 & 1 & 5 &(2,2)& 5 & 
11$\cdot$79 & 1 & 79 &(2,78)& 79 \\
 7$\cdot$127 & $3^{4}\cdot7^{2}$ & - &(6,42)& 7 &
13$\cdot$71 & $3^{3}$ & 61 &(12,10)& 61 \\ 
 5$\cdot$197 & $3^{3}\cdot11$ & - &(2,2)& 3 &
3$\cdot$331 & 1 & 3 & (2,6) &$3^{2}$ $(\ast)$ \\
3$\cdot$367 & 1 & 3 &(2,6)& 3 & 
17$\cdot$67 & $3^{7}\cdot11^{7}$ & 89 &(8,22),& $11\cdot89 \ (\ast)$ \\
 7$\cdot$163 & 1 & 19 & (6,18)& 19 &
  &  & &(8,22)&   \\ 
 19$\cdot$61 & $3^{3}\cdot7$ & 73 & (18,12) &73 &
 17$\cdot$71 & $3^{2}$ & 17 &(16,2)& 17 \\ 
 7$\cdot$173 & 1 & 7 &(6,2)& 7 &
 17$\cdot$73 & $3^{4}\cdot7\cdot37\cdot109$ & 5 &(4,4)& 5 \\ 
11$\cdot$113 & $5\cdot37$ & 41&(10,8) & 41 &
 3$\cdot$419 & 1 & 3 &(2,2) & 3 \\ 
 13$\cdot$97 & $7^{3}$ & 5 &(4,4) & $5\cdot7^{2}\cdot97$ &
 31$\cdot$41 & $3^{3}\cdot5^{6}$ & 7& (6,2)&$7\cdot$ \\ 
  & & & (6,6) & &      & &11 & (10,10) & $11\cdot$\\
  & & 97 & (12,96) & &      & &31 & (30,10) & $31$\\
 13$\cdot$101 & $3\cdot5^2$ & 31 &(6,10) & 31 &
17$\cdot$79 & $5$ & 17 &(16,2)& 17 \\ 
5$\cdot$271 & $3^{3}\cdot5$ & 37&(4,18) & 37 &
 5$\cdot$277 & $3^{2}\cdot139$ & 5 \ 7 &(4,4) (2,6)& $5\cdot7$ \\ 
19$\cdot$73 & $3^{4}\cdot7\cdot101$ & 17 & (18,8)& $17^{2}\cdot19\cdot37$ &
 7$\cdot$199 & 1 & 5&(2,2) & 5 \\ 
  & & 19 & (18,18)& & & & & & \\
  & & 37 & (18,36)& & & & & & \\
 5$\cdot$293 & $3^{2}\cdot7^{2}$ & - & (4,4)& $3^{2}$ &
 7$\cdot$211 & $3^{2}\cdot5^{2}\cdot$ & 11 &(2,10) (6,42)& $7 \cdot 11$ \\ 
 & & & & & & $7^{2}$& & & \\
3$\cdot$503 & 1 & 3 & (2,2) & 3 &
 17$\cdot$89 & $11^{3}\cdot17\cdot41$ & 13 &(4,4) (16,8)& $13\cdot17$ \\ 
37$\cdot$43 & $3^{26}\cdot7^{8}\cdot11\cdot$ & 43 &(6,42)& 43 &
 3$\cdot$541 & 1 & 13 &(2,12)& 13 \\ 
 & $19\cdot487$& & & & & & & & \\
 3$\cdot$547 & 1 & 5 &(2,2)& 5 &
 13$\cdot$127 & $3^{3}\cdot7$ & 5 & (12,2)& $5^{2}$ \\ 
7$\cdot$241 & 1 & 13 &(6,12)& 13 &
 5$\cdot$347 & $3\cdot29$ & 5 &(4,2) & 5 \\ 
37$\cdot$47 & $23^{5}$ & 5 &(4,2)& 5 &
 17$\cdot$103 & $3^{7}\cdot17^{7}$ & - & (8,6) (16,34) &$3^{2}\cdot17 \ (\ast)$ \\ 
3$\cdot$587 & 1 & 7 &(2,2)& 7 &
 5$\cdot$353 & $3^{2}\cdot59$ & - &(2,2)& 3  \\ 
5$\cdot$373 & $3^{2}\cdot11\cdot17$ & 5 &(4,4) & 5  &
 11$\cdot$173 & 3 & 173 & (2,172)&173 \\ 
17$\cdot$113 & $3^{3}\cdot19$ & 17 & (16,16)&$17^{3}\cdot29$ &
 13$\cdot$149 & $3^{2}\cdot5^{2}\cdot7$ & 109&(12,4) & $3\cdot109$ \\
 & & 29 & (4,28) & & & & & (6,2) & \\
\hline \hline
\end{tabular}
\end{center}

\begin{center} A{\scriptsize CKNOWLEDGMENT} \end{center}
The author wishes to thank Professor Lawrence Washington for his advice and guidance during the writing of the Thesis that this paper is based on, as well as the Cy-Tera of the Cyprus Institute and the Mathematics Department of the University of Maryland College Park, for their permission to run parts of this project on their machines.

\vspace{10mm}

\begin{flushleft}Eleni Agathocleous \\ email: agathocleous.e@unic.ac.cy \end{flushleft}

\end{document}